\documentclass{amsart}
\usepackage{amsmath,amssymb}
\usepackage{yhmath}

\newcommand{\bdis}{\begin{displaymath}}
\newcommand{\edis}{\end{displaymath}}
\newcommand{\be}{\begin{equation}}
\newcommand{\ee}{\end{equation}}
\newcommand{\mbb}{\mathbb}
\newcommand{\mcal}{\mathcal}

\newcommand{\vp}{\varphi}

\newcommand{\zf}{\zeta\left(\frac{1}{2}+it\right)}
 
\newcommand{\zfn}{\zeta\left(\frac{1}{2}+it_\nu\right)}   
\newcommand{\zfnn}{\zeta\left(\frac{1}{2}+it_{\nu+1}\right)}

\newcommand{\FR}{\frac{x^n+y^n}{z^n}}

\DeclareMathOperator{\im}{Im}

\theoremstyle{definition}

\theoremstyle{remark}
\newtheorem{remark}[]{Remark}

\newtheorem*{mydef11}{{\bf Theorem 1}}

\newtheorem*{mydef12}{{\bf Theorem 2}}

\newtheorem*{mydef13}{{\bf Theorem 3}}

\newtheorem*{mydef14}{{\bf Theorem 4}}

\newtheorem*{mydef51}{{\bf Lemma 1}}

\newtheorem*{mydef52}{{\bf Lemma 2}}

\newtheorem*{mydef53}{{\bf Lemma 3}}

\newtheorem*{mydef54}{{\bf Lemma 4}}

\newtheorem*{mydef81}{{\bf Property 1}}

\newtheorem*{mydef82}{{\bf Property 2}}

\numberwithin{equation}{section}

\begin{document}

\title[Jacob's ladders, E. C. Titchmarsh's hypothesis  \dots]{Jacob's ladders, E. C. Titchmarsh's hypothesis (1934) and new $\zeta$-equivalents of the Fermat-Wiles theorem or connections between Fermat's rationals and the Gram's sequence}

\author{Jan Moser}

\address{Department of Mathematical Analysis and Numerical Mathematics, Comenius University, Mlynska Dolina M105, 842 48 Bratislava, SLOVAKIA}

\email{jan.mozer@fmph.uniba.sk}


\begin{abstract}
In connection of our proof (1980) of the Titchmarsh's hypothesis (1934), we have obtained two asymptotic formulae (1991). In this paper we obtain three new $\zeta$-equivalents of the Fermat-Wiles theorem based on the mentioned asymptotic formulae. 
\end{abstract}
\maketitle

\section{Introduction} 

This paper is a continuation of the series of our ten papers about the points of contact between the Riemann's zeta-function and the Fermat-Wiles theorem.  

\subsection{}  

Let us remind the definition of the Riemann's zeta function 

\be \label{1.1} 
\zeta(s)=\prod_{p}\frac{1}{1-\frac{1}{p^s}},\ s=\sigma+it,\ \sigma>1, 
\ee 
where $p$ runs over the set of all primes, and analytic continuation of this function on all complex $s$ except $s=1$. Of course, this definition is based on the Euler's fundamental identity 
\be \label{1.2} 
\prod_{p}\frac{1}{1-\frac{1}{p^x}}=\sum_{n=1}^\infty \frac{1}{n^x},\ x>1. 
\ee  
Riemann defined also the following real-valued function 
\be \label{1.3} 
Z(t)=e^{i\vartheta(t)}\zf
\ee 
where\footnote{See \cite{14}, (35), (44), (62); \cite{16}, p. 98.}
\be \label{1.4} 
\begin{split} 
& \vartheta(t)=-\frac{t}{2}\ln\pi+\im\Gamma\left(\frac 14+i\frac{t}{2}\right)= \\ 
& \frac{t}{2}\ln\frac{t}{2\pi}-\frac{t}{2}-\frac{\pi}{8}+\mcal{O}\left(\frac 1t\right). 
\end{split} 
\ee 

\begin{remark}
It follows from (\ref{1.3}), that the properties of the basic signal (\ref{1.3}) generated by the function $\zf$ are connected with the distribution of the primes. This fact is to be regarded as a pleasant circumstance from the point of view of the Pythagorean philosophy of the Universe. 
\end{remark} 

\subsection{} 

Let us remind the following facts: 
\begin{itemize}
	\item[(A)] Titchmarsh presented in 1934 the hypothesis\footnote{\cite{16}, p. 105.}: there is $A>0$ such that 
	\be \label{1.5} 
	\sum_{\nu=M+1}^N Z^2(t_\nu)Z^2(t_{\nu+1})=\mcal{O}(N\ln^AN), 
	\ee 
	where, if $t_N=T$, 
	\be \label{1.6} 
	N\sim \frac{1}{2\pi}T\ln T, 
	\ee 
	and $M$ is sufficiently big fixed number, and next, the sequence $\{ t_\nu\}$ (the Gram's sequence) is defined by the condition\footnote{Comp. \cite{16}, p. 99.} 
	\be \label{1.7} 
	\vartheta(t_\nu)=\pi\nu,\ \nu = 1,2,\dots
	\ee 
\end{itemize} 

\begin{remark}
Namely, E. C. Titchmarsh finished his paper with this sentence: "\dots but there are additional complications, and the conjecture not been verified."
\end{remark} 

\begin{itemize}
	\item[(B)] We have proved this Titchmarsh's hypothesis in 1980 with $A=4$, see \cite{3}. 
	\item[(C)] Next, in 1991, we have derived the following asymptotic formula 
	\be \label{1.8} 
	\sum_{T\leq t_\nu\leq 2T} Z^2(t_\nu)Z^2(t_{\nu+1})=\frac{3}{4\pi^5}T\ln^5T+\mcal{O}(T\ln^4T) 
	\ee 
	(see \cite{4}, (2.4), (2.6) and (2.10)). 
\end{itemize} 

\begin{remark}
That means, we have obtained the final result: the order of the Titchmarsh's sum is 
\be \label{1.9} 
N\ln^4N, 
\ee 
comp. (\ref{1.6}) and (\ref{1.8}). 
\end{remark} 

\subsection{} 

In this paper we obtain, for example, the functional (cross-bred) based on our asymptotic formula (\ref{1.8}). Namely, it is true that: 
\be \label{1.10} 
\begin{split}
& \lim_{\tau\to\infty}\frac{1}{\tau}
\left\{
\int_{\frac{4\pi^5}{3\zeta^5(2\sigma)}x\tau}^{[\frac{4\pi^5}{3\zeta^5(2\sigma)}x\tau]^1}|\zeta(\sigma+it)|^2{\rm d}t
\right\}^5 \times \\ 
& \left\{
\int_{\frac{4\pi^5}{3\zeta^5(2\sigma)}x\tau}^{[\frac{4\pi^5}{3\zeta^5(2\sigma)}x\tau]^1}\left|\zf\right|^2{\rm d}t
\right\}^{-5} \times \\ 
& \sum_{t_\nu\geq \frac{4\pi^5}{3\zeta^5(2\sigma)}x\tau}^{\leq \frac{8\pi^5}{3\zeta^5(2\sigma)}x\tau}\left|\zfn\right|^2\left|\zfnn\right|^2 = x,\ x>0
\end{split}
\ee 
for every fixed $\sigma\geq \frac 12+\epsilon$ and every fixed small $\epsilon>0$. 

\subsection{} 

Next, in the special case of the Fermat's rationals 
\be \label{1.11} 
x\to \FR,\ x,y,z,n\in\mbb{N},\ n\geq 3, 
\ee 
it follows from the equation (\ref{1.10}) that $\zeta$-condition 

\be \label{1.12} 
\begin{split}
& \lim_{\tau\to\infty}\frac{1}{\tau}
\left\{
\int_{\frac{4\pi^5}{3\zeta^5(2\sigma)}\FR\tau}^{[\frac{4\pi^5}{3\zeta^5(2\sigma)}\FR\tau]^1}|\zeta(\sigma+it)|^2{\rm d}t
\right\}^5 \times \\ 
& \left\{
\int_{\frac{4\pi^5}{3\zeta^5(2\sigma)}\FR\tau}^{[\frac{4\pi^5}{3\zeta^5(2\sigma)}\FR\tau]^1}\left|\zf\right|^2{\rm d}t
\right\}^{-5} \times \\ 
& \sum_{t_\nu\geq \frac{4\pi^5}{3\zeta^5(2\sigma)}\FR\tau}^{\leq \frac{8\pi^5}{3\zeta^5(2\sigma)}\FR\tau}\left|\zfn\right|^2\left|\zfnn\right|^2 \not=1
\end{split}
\ee 
on the class of all Fermat's rationals (\ref{1.11}) represents new $\zeta$-equivalent of the Fermat-Wiles theorem for every fixed $\sigma\geq \frac 12+\epsilon$ and every small fixed $\epsilon>0$. 

\subsection{} 

Further, in this paper we prove the following result 
\be \label{1.13} 
\frac{\int_{\overset{r-1}{T}}^{\overset{r}{T}}|\zf|^2{\rm d}t}{\int_{\overset{r-1}{T}}^{\overset{r}{T}}|S_1(t)|^{2l}{\rm d}t}=\frac{1}{\bar{c}(l)}\ln\overset{r-1}{T}+\mcal{O}(1),\ r=1,\dots,k,\ T\to\infty
\ee 
for every fixed $k,l\in\mbb{N}$, where 
\be \label{1.14} 
S_1(t)=\frac{1}{\pi}\int_0^t\arg\left\{\zeta\left(\frac 12+iu\right)\right\}{\rm d}u. 
\ee 

Next, on the basis of the formula (\ref{1.13}), we obtain the following functional (cross-bred) 
\be \label{1.15} 
\begin{split}
& \lim_{\tau\to\infty}\frac{1}{\tau}
\left\{
\int_{\frac{4\pi^5}{3\{\bar{c}(l)\}^5}x\tau}^{[\frac{4\pi^5}{3\{\bar{c}(l)\}^5}x\tau]^1}|S_1(t)|^{2l}{\rm d}t
\right\}^5 \times \\ 
& \left\{
\int_{\frac{4\pi^5}{3\{\bar{c}(l)\}^5}x\tau}^{[\frac{4\pi^5}{3\{\bar{c}(l)\}^5}x\tau]^1}\left|\zf\right|^2{\rm d}t
\right\}^{-5} \times \\ 
& \sum_{t_\nu\geq \frac{4\pi^5}{3\{\bar{c}(l)\}^5}x\tau}^{\leq \frac{8\pi^5}{3\{\bar{c}(l)\}^5}x\tau}\left|\zfn\right|^2\left|\zfnn\right|^2 = x,\ x>0 
\end{split}
\ee 
for every fixed $l\in\mbb{N}$. 

\begin{remark}
Of course, next $\zeta$-equivalent of the Fermat-Wiles theorem follows from (\ref{1.15}) as an analogue of the eq. (\ref{1.12}). 
\end{remark}

\section{Jacob's ladders: notions and basic geometrical properties}  

\subsection{} 

In this paper we use the following notions of our works \cite{5} -- \cite{9}: 
\begin{itemize}
\item[{\tt (a)}] Jacob's ladder $\vp_1(T)$, 
\item[{\tt (b)}] direct iterations of Jacob's ladders 
\bdis 
\begin{split}
	& \vp_1^0(t)=t,\ \vp_1^1(t)=\vp_1(t),\ \vp_1^2(t)=\vp_1(\vp_1(t)),\dots , \\ 
	& \vp_1^k(t)=\vp_1(\vp_1^{k-1}(t))
\end{split}
\edis 
for every fixed natural number $k$, 
\item[{\tt (c)}] reverse iterations of Jacob's ladders 
\be \label{2.1}  
\begin{split}
	& \vp_1^{-1}(T)=\overset{1}{T},\ \vp_1^{-2}(T)=\vp_1^{-1}(\overset{1}{T})=\overset{2}{T},\dots, \\ 
	& \vp_1^{-r}(T)=\vp_1^{-1}(\overset{r-1}{T})=\overset{r}{T},\ r=1,\dots,k, 
\end{split} 
\ee   
where, for example, 
\be \label{2.2} 
\vp_1(\overset{r}{T})=\overset{r-1}{T}
\ee  
for every fixed $k\in\mbb{N}$ and every sufficiently big $T>0$. We also use the properties of the reverse iterations listed below.  
\be \label{2.3}
\overset{r}{T}-\overset{r-1}{T}\sim(1-c)\pi(\overset{r}{T});\ \pi(\overset{r}{T})\sim\frac{\overset{r}{T}}{\ln \overset{r}{T}},\ r=1,\dots,k,\ T\to\infty,  
\ee 
\be \label{2.4} 
\overset{0}{T}=T<\overset{1}{T}(T)<\overset{2}{T}(T)<\dots<\overset{k}{T}(T), 
\ee 
and 
\be \label{2.5} 
T\sim \overset{1}{T}\sim \overset{2}{T}\sim \dots\sim \overset{k}{T},\ T\to\infty.   
\ee  
\end{itemize} 

\begin{remark}
	The asymptotic behaviour of the points 
	\bdis 
	\{T,\overset{1}{T},\dots,\overset{k}{T}\}
	\edis  
	is as follows: at $T\to\infty$ these points recede unboundedly each from other and all together are receding to infinity. Hence, the set of these points behaves at $T\to\infty$ as one-dimensional Friedmann-Hubble expanding Universe. 
\end{remark}  

\subsection{} 

Let us remind that we have proved\footnote{See \cite{9}, (3.4).} the existence of almost linear increments 
\be \label{2.6} 
\begin{split}
& \int_{\overset{r-1}{T}}^{\overset{r}{T}}\left|\zf\right|^2{\rm d}t\sim (1-c)\overset{r-1}{T}, \\ 
& r=1,\dots,k,\ T\to\infty,\ \overset{r}{T}=\overset{r}{T}(T)=\vp_1^{-r}(T)
\end{split} 
\ee 
for the Hardy-Littlewood integral (1918) 
\be \label{2.7} 
J(T)=\int_0^T\left|\zf\right|^2{\rm d}t. 
\ee  

For completeness, we give here some basic geometrical properties related to Jacob's ladders. These are generated by the sequence 
\be \label{2.8} 
T\to \left\{\overset{r}{T}(T)\right\}_{r=1}^k
\ee 
of reverse iterations of the Jacob's ladders for every sufficiently big $T>0$ and every fixed $k\in\mbb{N}$. 

\begin{mydef81}
The sequence (\ref{2.8}) defines a partition of the segment $[T,\overset{k}{T}]$ as follows 
\be \label{2.9} 
|[T,\overset{k}{T}]|=\sum_{r=1}^k|[\overset{r-1}{T},\overset{r}{T}]|
\ee 
on the asymptotically equidistant parts 
\be \label{2.10} 
\begin{split}
& \overset{r}{T}-\overset{r-1}{T}\sim \overset{r+1}{T}-\overset{r}{T}, \\ 
& r=1,\dots,k-1,\ T\to\infty. 
\end{split}
\ee 
\end{mydef81} 

\begin{mydef82}
Simultaneously with the Property 1, the sequence (\ref{2.8}) defines the partition of the integral 
\be \label{2.11} 
\int_T^{\overset{k}{T}}\left|\zf\right|^2{\rm d}t
\ee 
into the parts 
\be \label{2.12} 
\int_T^{\overset{k}{T}}\left|\zf\right|^2{\rm d}t=\sum_{r=1}^k\int_{\overset{r-1}{T}}^{\overset{r}{T}}\left|\zf\right|^2{\rm d}t, 
\ee 
that are asymptotically equal 
\be \label{2.13} 
\int_{\overset{r-1}{T}}^{\overset{r}{T}}\left|\zf\right|^2{\rm d}t\sim \int_{\overset{r}{T}}^{\overset{r+1}{T}}\left|\zf\right|^2{\rm d}t,\ T\to\infty. 
\ee 
\end{mydef82} 

It is clear, that (\ref{2.10}) follows from (\ref{2.3}) and (\ref{2.5}) since 
\be \label{2.14} 
\overset{r}{T}-\overset{r-1}{T}\sim (1-c)\frac{\overset{r}{T}}{\ln \overset{r}{T}}\sim (1-c)\frac{T}{\ln T},\ r=1,\dots,k, 
\ee  
while our eq. (\ref{2.13}) follows from (\ref{2.6}) and (\ref{2.5}). 

\section{The equivalent of the Fermat-Wiles theorem generated by our first asymptotic formula (1991)}

\subsection{} 

We use the following $\zeta$-formulas: 

\begin{itemize}
	\item[(A)] The quotient formula\footnote{See \cite{11}, (3.11).} 
	\be \label{3.1} 
	\begin{split} 
	& \frac{\int_{\overset{r-1}{T}}^{\overset{r}{T}}|\zf|^2{\rm d}t}{\int_{\overset{r-1}{T}}^{\overset{r}{T}}|\zeta(\sigma+it)|^{2}{\rm d}t}=\frac{1}{
	\zeta(2\sigma)}\ln\overset{r-1}{T}+\mcal{O}(1), \\ 
	& r=1,\dots,k,\ T\to\infty 
	\end{split}  
	\ee 
	for every fixed $k\in\mbb{N}$, $\sigma\geq \frac 12+\epsilon$, that is based on the Hardy-Littlewood formula (1922), comp. \cite{11}, (3.1), and our almost linear formula, comp. \cite{9}, (3.6). 
	\item[(B)] The asymptotic formula (\ref{1.8}), (\ref{1.9}) 
	\be \label{3.2} 
	\begin{split} 
	& \sum_{t_\nu\geq T}^{t_\nu\leq 2T}\left|\zfn\right|^2\left|\zfnn\right|^2=\\ 
	& \frac{3}{4\pi^5}T\ln^5T\left\{1+\mcal{O}\left(\frac{1}{\ln T}\right)\right\}
	\end{split}
	\ee  
	for the Titchmarsh's sum. 
	\end{itemize} 
	
Since (see (\ref{3.1}), $r=1$) 
\be \label{3.3} 
\ln T=\zeta(2\sigma)\frac{\int_{T}^{\overset{1}{T}}|\zf|^2{\rm d}t}{\int_{T}^{\overset{1}{T}}|\zeta(\sigma+it)|^{2}{\rm d}t}\left\{1+\mcal{O}\left(\frac{1}{\ln T}\right)\right\}, 
\ee 
then it follows from (\ref{3.2}) that 
\be \label{3.4} 
\begin{split}
& \left\{\frac{\int_{T}^{\overset{1}{T}}|\zf|^2{\rm d}t}{\int_{T}^{\overset{1}{T}}|\zeta(\sigma+it)|^{2}{\rm d}t}\right\}^5 \times \\ 
& \sum_{t_\nu\geq T}^{t_\nu\leq 2T}\left|\zfn\right|^2\left|\zfnn\right|^2= \\ 
& \frac{3\zeta^5(2\sigma)}{4\pi^5}T\left\{1+\mcal{O}\left(\frac{1}{\ln T}\right)\right\}. 
\end{split}
\ee 

Now, the substitution 
\be \label{3.5} 
T=\frac{4\pi^5}{3\zeta^5(2\sigma)}x\tau ,\ x>0 
\ee 
into the eq. (\ref{3.4}) gives the following functional (as the cross-breed of our asymptotic formula (\ref{1.3})). 

\begin{mydef51}
\be \label{3.6} 
\begin{split}
& \lim_{\tau\to\infty}\frac{1}{\tau}\left\{
\int_{\frac{4\pi^5}{3\zeta^5(2\sigma)}x\tau}^{[\frac{4\pi^5}{3\zeta^5(2\sigma)}x\tau]^1}|\zeta(\sigma+it)|^2{\rm d}t
\right\}^5 \times \\ 
& \left\{
\int_{\frac{4\pi^5}{3\zeta^5(2\sigma)}x\tau}^{[\frac{4\pi^5}{3\zeta^5(2\sigma)}x\tau]^1}\left|\zf\right|^2{\rm d}t
\right\}^{-5} \times \\ 
& \sum_{t_\nu\geq \frac{4\pi^5}{3\zeta^5(2\sigma)}x\tau}^{\leq \frac{8\pi^5}{3\zeta^5(2\sigma)}x\tau}\left|\zfn\right|^2\left|\zfnn\right|^2 = x 
\end{split}
\ee 
for every fixed 
\bdis 
x>0,\ \sigma\geq\frac 12+\epsilon. 
\edis 
\end{mydef51} 

\subsection{} 

Now, in the special case\footnote{See (\ref{1.11}).} 
\bdis 
x\to\FR 
\edis  wee obtain the following Lemma from (\ref{3.6}). 

\begin{mydef52}
\be \label{3.7} 
\begin{split}
& \lim_{\tau\to\infty}\frac{1}{\tau}
\left\{
\int_{\frac{4\pi^5}{3\zeta^5(2\sigma)}\FR\tau}^{[\frac{4\pi^5}{3\zeta^5(2\sigma)}\FR\tau]^1}|\zeta(\sigma+it)|^2{\rm d}t
\right\}^5 \times \\ 
& \left\{
\int_{\frac{4\pi^5}{3\zeta^5(2\sigma)}\FR\tau}^{[\frac{4\pi^5}{3\zeta^5(2\sigma)}\FR\tau]^1}\left|\zf\right|^2{\rm d}t
\right\}^{-5} \times \\ 
& \sum_{t_\nu\geq \frac{4\pi^5}{3\zeta^5(2\sigma)}\FR\tau}^{\leq \frac{8\pi^5}{3\zeta^5(2\sigma)}\FR\tau}\left|\zfn\right|^2\left|\zfnn\right|^2 = \FR
\end{split}
\ee 
for every fixed Fermat's rational and every fixed $\sigma\geq \frac 12+\epsilon$. 
\end{mydef52} 

Consequently, we have the following result 

\begin{mydef11}
The $\zeta$-condition 
\be \label{3.8} 
\begin{split}
& \lim_{\tau\to\infty}\frac{1}{\tau}
\left\{
\int_{\frac{4\pi^5}{3\zeta^5(2\sigma)}\FR\tau}^{[\frac{4\pi^5}{3\zeta^5(2\sigma)}\FR\tau]^1}|\zeta(\sigma+it)|^2{\rm d}t
\right\}^5 \times \\ 
& \left\{
\int_{\frac{4\pi^5}{3\zeta^5(2\sigma)}\FR\tau}^{[\frac{4\pi^5}{3\zeta^5(2\sigma)}\FR\tau]^1}\left|\zf\right|^2{\rm d}t
\right\}^{-5} \times \\ 
& \sum_{t_\nu\geq \frac{4\pi^5}{3\zeta^5(2\sigma)}\FR\tau}^{\leq \frac{8\pi^5}{3\zeta^5(2\sigma)}\FR\tau}\left|\zfn\right|^2\left|\zfnn\right|^2 \not=1 
\end{split}
\ee 
on the class of all Fermat's rationals represents the next $\zeta$-equivalent of the Fermat-Wiles theorem for every fixed $\sigma\geq \frac 12+\epsilon$. 
\end{mydef11} 

\section{The second equivalent of the Fermat-Wiles theorem generated by our first asymptotic formula (1991)} 

\subsection{} 

Let us remind the following Selberg's formula, see \cite{15}, p. 130: 

\be \label{4.1} 
\begin{split}
& \int_T^{T+H}|S_1(t)|^{2l}{\rm d}t=\bar{c}(l)H+\mcal{O}\left(\frac{H}{\ln T}\right), \\ 
& T^a\leq H\leq T,\ \frac 12< a\leq 1,\ \bar{c}(l)>0,\ l\in \mbb{N}. 
\end{split}
\ee 
We use in this section: 

\begin{itemize}
	\item[(A)] The formula (\ref{4.1}) in the case 
	\be \label{4.2} 
	T=\overset{r-1}{T},\ H=\overset{r}{T}-\overset{r-1}{T}, 
	\ee  
	that gives\footnote{See (\ref{2.5}), (\ref{2.4}).} 
	\be \label{4.3} 
	\begin{split}
	& \int_{\overset{r-1}{T}}^{\overset{r}{T}}|S_1(t)|^{2l}{\rm d}t=\bar{c}(l)(\overset{r}{T}-\overset{r-1}{T})+\mcal{O}\left(\frac{T}{\ln^2 T}\right), \\ 
	& r=1,\dots,k,\ T\to\infty
	\end{split}
	\ee 
	for every fixed $k,l\in\mbb{N}$. 
	\item[(B)] Our almost linear formula, see \cite{9}, (3.4) 
	\be \label{4.4} 
	\begin{split}
	& \int_{\overset{r-1}{T}}^{\overset{r}{T}}\left|\zf\right|^2{\rm d}t=(1-c)\overset{r-1}{T}+\mcal{O}(T^{1/3+\delta}), \\ 
	& r=1,\dots,k,\ T\to\infty, 
	\end{split}
	\ee 
	where $\delta>0$ being sufficiently small. 
\end{itemize}

Consequently, we obtain the next theorem for the quotient of formulae (\ref{4.4}) and (\ref{4.3}) in the exactly same way as we did in our paper \cite{11}, (3.7) -- (3.10). 

\begin{mydef12}
It is true, that  
\be\label{4.5} 
\begin{split}
& \frac{\int_{\overset{r-1}{T}}^{\overset{r}{T}}\left|\zf\right|^2{\rm d}t}{\int_{\overset{r-1}{T}}^{\overset{r}{T}}|S_1(t)|^{2l}{\rm d}t}=\frac{1}{\bar{c}(l)}\ln \overset{r-1}{T}+\mcal{O}(1), \\ 
& r=1,\dots,k,\ T\to\infty, 
\end{split}
\ee 
for every fixed $k,l\in\mbb{N}$. 
\end{mydef12} 

\subsection{} 

Here, we use the formula (\ref{4.5}) with $r=1$ in the form\footnote{Comp. (\ref{3.3}).} 
\be \label{4.6} 
\ln T=\bar{c}(l)\frac{\int_{T}^{\overset{1}{T}}\left|\zf\right|^2{\rm d}t}{\int_{T}^{\overset{1}{T}}|S_1(t)|^{2l}{\rm d}t}\left\{1+\mcal{O}\left(\frac{1}{\ln T}\right)\right\}. 
\ee 
Now, the cross-breeding of (\ref{3.2}) and (\ref{4.6}) gives us the following 
\be \label{4.7} 
\begin{split}
& \left\{\int_{T}^{\overset{1}{T}}|S_1(t)|^{2l}{\rm d}t\right\}^5\times \left\{\int_{T}^{\overset{1}{T}}\left|\zf\right|^2{\rm d}t\right\}^{-5}\times \\ 
& \sum_{t_\nu\geq T}^{t_\nu\leq 2T}\left|\zfn\right|^2\left|\zfnn\right|^2= \\ 
& \frac{3\{\bar{c}(l)\}^5}{4\pi^5}T\left\{1+\mcal{O}\left(\frac{1}{\ln T}\right)\right\} 
\end{split}
\ee 
for every fixed $l\in\mbb{N}$. 

Consequently, the substitution 
\be \label{4.8} 
T=\frac{4\pi^5}{3\{\bar{c}(l)\}^5}x\tau,\ x>0 
\ee  
into the eq. (\ref{4.7}) gives the following functional (as the cross-breed of our asymptotic formula (\ref{3.2}) for the Titchmarsh's sum). 

\begin{mydef53}
\be \label{4.9} 
\begin{split}
& \lim_{\tau\to\infty}\frac{1}{\tau}\left\{
\int_{\frac{4\pi^5}{3\{\bar{c}(l)\}^5}x\tau}^{[\frac{4\pi^5}{3\{\bar{c}(l)\}^5}x\tau]^1}|S_1(t)|^{2l}{\rm d}t
\right\}^5 \times \\ 
& \left\{
\int_{\frac{4\pi^5}{3\{\bar{c}(l)\}^5}x\tau}^{[\frac{4\pi^5}{3\{\bar{c}(l)\}^5}x\tau]^1}\left|\zf\right|^2{\rm d}t
\right\}^{-5} \times \\ 
& \sum_{t_\nu\geq \frac{4\pi^5}{3\{\bar{c}(l)\}^5}x\tau}^{\leq \frac{8\pi^5}{3\{\bar{c}(l)\}^5}x\tau}\left|\zfn\right|^2\left|\zfnn\right|^2 =x  
\end{split}
\ee 
for every fixed $x>0$ and $l\in\mbb{N}$. 
\end{mydef53} 

Finally, we obtain from (\ref{4.9}), exactly as we did it in (\ref{3.6}) -- (\ref{3.8}), the following result. 

\begin{mydef13}
The $\zeta$-condition 
\be \label{4.10} 
\begin{split}
& \lim_{\tau\to\infty}\frac{1}{\tau}\left\{
\int_{\frac{4\pi^5}{3\{\bar{c}(l)\}^5}\FR\tau}^{[\frac{4\pi^5}{3\{\bar{c}(l)\}^5}\FR\tau]^1}|S_1(t)|^{2l}{\rm d}t
\right\}^5 \times \\ 
& \left\{
\int_{\frac{4\pi^5}{3\{\bar{c}(l)\}^5}\FR\tau}^{[\frac{4\pi^5}{3\{\bar{c}(l)\}^5}\FR\tau]^1}\left|\zf\right|^2{\rm d}t
\right\}^{-5} \times \\ 
& \sum_{t_\nu\geq \frac{4\pi^5}{3\{\bar{c}(l)\}^5}\FR\tau}^{\leq \frac{8\pi^5}{3\{\bar{c}(l)\}^5}\FR\tau}\left|\zfn\right|^2\left|\zfnn\right|^2 \not=1
\end{split}
\ee 
on the class of all Fermat's rationals represents new $\zeta$-equivalent of the Fermat-Wiles theorem for every fixed $l\in\mbb{N}$. 
\end{mydef13} 

\section{The equivalent of the Fermat-Wiles theorem generated by our second asymptotic formula (1991)} 

\subsection{} 

Let us remind our second asymptotic formula\footnote{See \cite{4}, (2.4), (4.4).} 
\be \label{5.1} 
\sum_{t_\nu\geq T}^{t_\nu\leq 2T}\left|\zfn\right|^4=\frac{1}{4\pi^3}T\ln^5T+\mcal{O}(T\ln^4T). 
\ee  
We will use this one in the form 
\be \label{5.2} 
\sum_{t_\nu\geq T}^{t_\nu\leq 2T}\left|\zfn\right|^4=\frac{1}{4\pi^3}T\ln^5T\left\{1+\mcal{O}\left(\frac{1}{\ln T}\right)\right\}. 
\ee  
Next, the cross-breeding of this formula with the formula (\ref{3.3}) gives us the result 
\be \label{5.3} 
\begin{split}
& \left\{\frac
{\int_T^{\overset{1}{T}}|\zeta(\sigma+it)|^2{\rm d}t}
{\int_T^{\overset{1}{T}}|\zf|^2{\rm d}t}\right\}^5\times 
\sum_{t_\nu\geq T}^{t_\nu\leq 2T}\left|\zfn\right|^4= \\ 
& \frac{\zeta^5(2\sigma)}{4\pi^3}T\left\{1+\mcal{O}\left(\frac{1}{\ln T}\right)\right\} 
\end{split}
\ee 
for every fixed $\sigma\geq \frac 12+\epsilon$. 

Now, the substitution 
\be \label{5.4} 
T=\frac{4\pi^3}{\zeta^5(2\sigma)}x\tau,\ x>0 
\ee 
into the eq. (\ref{5.3}) gives to us the following functional (as the cross-breed of our formula (\ref{5.1})). 

\begin{mydef54}
\be \label{5.5} 
\begin{split}
& \lim_{\tau\to\infty}\frac{1}{\tau}\left\{\int_{\frac{4\pi^3}{\zeta^5(2\sigma)}x\tau}^{[\frac{4\pi^3}{\zeta^5(2\sigma)}x\tau]^1}|\zeta(\sigma+it)|^2{\rm d}t\right\}^5\times \\ 
& \left\{\int_{\frac{4\pi^3}{\zeta^5(2\sigma)}x\tau}^{[\frac{4\pi^3}{\zeta^5(2\sigma)}x\tau]^1}|\zf|^2{\rm d}t\right\}^{-5}\times \\ 
& \sum_{t_\nu\geq \frac{4\pi^3}{\zeta^5(2\sigma)}x\tau}^{t_\nu\leq \frac{8\pi^3}{\zeta^5(2\sigma)}x\tau}\left|\zfn\right|^4=x
\end{split}
\ee 
for every fixed 
\bdis 
x>0,\ \sigma\geq \frac 12+\epsilon. 
\edis 
\end{mydef54} 

Finally, we obtain from (\ref{5.5}) (exactly as in the case (\ref{3.6}) -- (\ref{3.8})) the following theorem. 

\begin{mydef14}
The $\zeta$-condition 
\be \label{5.6} 
\begin{split}
& \lim_{\tau\to\infty}\frac{1}{\tau}\left\{\int_{\frac{4\pi^3}{\zeta^5(2\sigma)}\FR\tau}^{[\frac{4\pi^3}{\zeta^5(2\sigma)}\FR\tau]^1}|\zeta(\sigma+it)|^2{\rm d}t\right\}^5\times \\ 
& \left\{\int_{\frac{4\pi^3}{\zeta^5(2\sigma)}\FR\tau}^{[\frac{4\pi^3}{\zeta^5(2\sigma)}\FR\tau]^1}|\zf|^2{\rm d}t\right\}^{-5}\times \\ 
& \sum_{t_\nu\geq \frac{4\pi^3}{\zeta^5(2\sigma)}\FR\tau}^{t_\nu\leq \frac{8\pi^3}{\zeta^5(2\sigma)}\FR\tau}\left|\zfn\right|^4\not= 1
\end{split}
\ee 
on the class of all Fermat's rationals represents the next $\zeta$-equivalent of the Fermat-Wiles theorem for every fixed $\sigma\geq \frac 12+\epsilon$. 
\end{mydef14} 

\section{Continuation of the chain of equivalences from \cite{13}, (6.7) by means of new functionals of this paper} 

\subsection{} 

Namely, the functionals (\ref{3.6}), (\ref{4.9}) and (\ref{5.5}) give the following continuation of the chain of equivalences mentioned above: 
\be \label{6.1} 
\begin{split}
& \sim \left\{\int_{\frac{4\pi^5}{3\zeta^5(2\sigma)}x\tau}^{[\frac{4\pi^5}{3\zeta^5(2\sigma)}x\tau]^1}|\zeta(\sigma+it)|^2{\rm d}t\right\}^5\times \left\{\int_{\frac{4\pi^5}{3\zeta^5(2\sigma)}x\tau}^{[\frac{4\pi^5}{3\zeta^5(2\sigma)}x\tau]^1}|\zf|^2{\rm d}t\right\}^{-5}\times \\ 
& \sum_{t_\nu\geq \frac{4\pi^5}{3\zeta^5(2\sigma)}x\tau}^{t_\nu\leq \frac{8\pi^5}{3\zeta^5(2\sigma)}x\tau}\left|\zfn\right|^2\left|\zfnn\right|^2\sim \\ 
& \left\{
\int_{\frac{4\pi^5}{3\{\bar{c}(l)\}^5}x\tau}^{[\frac{4\pi^5}{3\{\bar{c}(l)\}^5}x\tau]^1}|S_1(t)|^{2l}{\rm d}t
\right\}^5 \times  
\left\{
\int_{\frac{4\pi^5}{3\{\bar{c}(l)\}^5}x\tau}^{[\frac{4\pi^5}{3\{\bar{c}(l)\}^5}x\tau]^1}\left|\zf\right|^2{\rm d}t
\right\}^{-5} \times \\ 
& \sum_{t_\nu\geq \frac{4\pi^5}{3\{\bar{c}(l)\}^5}x\tau}^{t_\nu\leq \frac{8\pi^5}{3\{\bar{c}(l)\}^5}x\tau}\left|\zfn\right|^2\left|\zfnn\right|^2\sim \\ 
& \left\{\int_{\frac{4\pi^3}{\zeta^5(2\sigma)}x\tau}^{[\frac{4\pi^3}{\zeta^5(2\sigma)}x\tau]^1}|\zeta(\sigma+it)|^2{\rm d}t\right\}^5\times  
 \left\{\int_{\frac{4\pi^3}{\zeta^5(2\sigma)}x\tau}^{[\frac{4\pi^3}{\zeta^5(2\sigma)}x\tau]^1}|\zf|^2{\rm d}t\right\}^{-5}\times \\ 
& \sum_{t_\nu\geq \frac{4\pi^3}{\zeta^5(2\sigma)}x\tau}^{t_\nu\leq \frac{8\pi^3}{\zeta^5(2\sigma)}x\tau}\left|\zfn\right|^4
\end{split}
\ee 
for every fixed 
\bdis 
x>0,\ l\in\mbb{N}, \ \sigma\geq \frac 12+\epsilon. 
\edis 

\begin{remark}
Of course, the set of (finite) chains, that is composed as the union of two parts 
\bdis 
\{[13],(6.7)\},\ \{(6.1)\}
\edis 
represents, for example in the case 
\bdis 
x\in [1,\alpha],\ \alpha>1, 
\edis 
the continuum set of different (not equivalent) chains for every fixed 
\bdis 
l\in\mbb{N},\ \sigma\geq \frac 12+\epsilon,\ \alpha>1, 
\edis 
comp. \cite{13}, Corollary 1. 
\end{remark}

\subsection{} 

In connection with the Remark 1 about the Pythagorean philosophy of the Universe we give here the G.H. Hardy's expression, \cite{2}, the end of the $22^{th}$ part: \\  \emph{I believe that mathematical reality lies outside us, that our function is to discover or observe it, and that the theorems which we prove, and which we describe grandiloquently as our creations, are simply our notes of our observations}.

I would like to thank Michal Demetrian for his moral support of my study of Jacob's ladders.

\end{document}